\documentstyle[twoside]{paper}

\font\smallgothic=eufm10 at 9truept
\newcommand{\Gm}{\mbox{\gothic m}}
\newcommand{\Gp}{\mbox{\gothic p}}
\newcommand{\Gq}{\mbox{\gothic q}}
\newcommand{\GT}{\mbox{\gothic T}}

\newcommand{\Gpp}{\mbox{\smallgothic p}}

\begin{document}
\begin{center}
{\Huge Explicit representations

\vspace{12pt}
of biorthogonal polynomials}

\vspace{16pt}
{\Large Arieh Iserles\footnote{Department of Applied Mathematics and
Theoretical Physics, University of Cambridge, Cambridge, England.}}
{\large and} {\Large Syvert P. N\o rsett\footnote{Institute of Mathematical
Sciences, Norwegian Institute of Technology, Trondheim, Norway.}}

\vspace{28pt}
\parbox[t]{130mm}{{\bf Abstract}  Given a parametrised weight function
$\omega(x,\mu)$ such that the quotients of its consecutive  moments
are M\"obius maps, it is possible to express the underlying
biorthogonal polynomials in a closed form \cite{IN2}. In the present
paper we address ourselves to two related issues. Firstly, we
demonstrate that, subject to additional assumptions, every such $\omega$
obeys (in $x$) a linear differential equation whose solution is a
generalized hypergeometric function. Secondly, using a generalization
of standard divided differences, we present a new explicit
representation of the underlying orthogonal polynomials.}
\end{center}

\vfill
\noindent AMS (MOS) {\em Mathematics Subject Classification.\/}
Primary 42C05; Secondary 33C45.
\eject

\section{Introduction}

Let us consider a one-parametric family of weight functions
$\omega(x,\mu)\geq0$, $x\geq0$, $\mu_1<\mu<\mu_+$. We define a
{\em biorthogonal polynomial\/} $p_n\in{\Bbb P}_n\times (\mu_-,\mu_+)^n$,
$p_n\not\equiv0$, by means of the orthogonality relations
\begin{equation}\label{1.1}
\int_0^1 p_n(x;\mu_1,\mu_2,\ldots,\mu_n)\omega(x,\mu_\ell) \D
x=0,\qquad \ell=1,2,\ldots,n,
\end{equation}
where $\mu_1,\mu_2,\ldots,\mu_n$ are distinct. (Here -- and elsewhere in
this paper -- ${\Bbb P}_n$ denotes the set of all $n$th degree
polynomials.)

The theory of general biorthogonal polynomials has been described in 
detail in \cite{IN1} and further developments have been reported in
\cite{IN2,IN3}. It possesses several interesting applications, ranging
from numerical methods for ordinary differential equations
\cite{IN0,IN5} to approximation theory \cite{IS1} to location of zeros
of polynomial transformations and orthogonal expansions
\cite{IN2,INS,IS2}. An exposition of biorthogonal functions and their
applications features in Claude Brezinski's monograph \cite{Brezinski}.
Herewith we revisit these elements of the theory
that are germane to the work of the present paper.

Biorthogonal polynomials exist and are unique (up to a nonzero
multiplicative constant) if and only if
$$\det\left[\begin{array}{cccc}
\Gm_0(\mu_1) & \Gm_1(\mu_1) & \cdots & \Gm_{m-1}(\mu_1)\\
\Gm_0(\mu_2) & \Gm_1(\mu_2) & \cdots & \Gm_{m-1}(\mu_2)\\
\vdots & \vdots & & \vdots\\
\Gm_0(\mu_m) & \Gm_1(\mu_m) & \cdots & \Gm_{m-1}(\mu_m)
	    \end{array}\right]\neq0$$
for every distinct $\mu_1,\mu_2,\ldots,\mu_n\in(\mu_-,\mu_+)$, where
$$\Gm_\ell(\mu)=\int_0^\infty x^\ell \omega(x,\mu)\D x$$
is the $\ell$th moment of $\omega$, $\ell\in\Zp$ \cite{IN1}. The zeros of $p_n$
need not, in general, reside in $(0,\infty)$. However, all the zeros
are in the support if $\omega$ is {\em strictly sign consistent
(SSC)\/}, i.e.\ if
$$\det \left[\begin{array}{cccc}
\omega(x_1,\mu_1) & \omega(x_2,\mu_1) & \cdots & \omega(x_m,\mu_1)\\
\omega(x_1,\mu_2) & \omega(x_2,\mu_2) & \cdots & \omega(x_m,\mu_2)\\
\vdots & \vdots & & \vdots\\
\omega(x_1,\mu_m) & \omega(x_2,\mu_m) & \cdots & \omega(x_m,\mu_m)
       \end{array}\right]\neq0$$
for every $n\geq1$ and all monotone sequences $0<x_1<x_2<\cdots<x_n$,
$\mu_-<\mu_1<\mu_2<\cdots<\mu_n<\mu_+$ \cite{IN1}.

Suppose that $p_n(\,\cdot\,;\mu_1,\mu_2,\ldots,\mu_n)$ is explicitly
known for every distinct $\mu_1,\ldots,\mu_n\in(\mu_-,\mu_+)$ and that
this explicit form can be extended to all
$\mu_1,\mu_2,\ldots,\mu_n\in\RR$.\footnote{It is immaterial whether
the definition of $\omega$ can be continued to all $\mu\in\RR$.} In
other words, given an arbitrary $n$-tuple $[
\mu_1,\mu_2,\ldots,\mu_n]\in\RR^d$ we can produce a $p_n$ or,
rephrasing again, we explicitly know a transformation $\GT:{\Bbb P}_n
\rightarrow {\Bbb P}_n$ such that 
\begin{equation}\label{1.2}
\GT \left\{ \prod_{\ell=1}^n (\,\cdot\,-\mu_\ell)\right\}
=p_n(\,\cdot\,;\mu_1,\mu_2,\ldots,\mu_n).
\end{equation}
If, in addition, we know that $\omega$ is SSC -- hence that all zeros
of $p_n$ reside in $(0,\infty)$ -- we deduce that $\GT$ maps
polynomials with all their zeros in $(\mu_-,\mu_+)$ into polynomials
with all their zeros in $(0,\infty)$ \cite{IN2}.

To date, the most efficacious mechanism for the generation of an explicit
form of $p_n$ has been presented in \cite{IN2}, leading to a long list
of interesting transformations of the form \R{1.2}. Thus, having
normalized $\Gm_0(\mu)\equiv1$, we asume that there exist real
sequences $\{\alpha_n\}_{n\in\Zp}$, $\{\beta_n\}_{n\in\Zp}$,
$\{\gamma_n\}_{n\in\Zp}$ and $\{\delta_n\}_{n\in\Zp}$ such that
\begin{equation}\label{1.3}
\alpha_k\delta_\ell-\beta_k\delta_\ell\neq0,\qquad
k=0,1,\ldots,\ell,\quad \ell\in\Zp,
\end{equation}
and
$$\Gm_n(\mu)=\prod_{\ell=0}^{n-1} \frac{\alpha_\ell+\mu\beta_\ell}
{\gamma_\ell+\mu \delta_\ell},\qquad n\in\Zp.$$
In other words,
\begin{equation}\label{1.4}
\Gm_{n+1}(\mu)=\frac{\alpha_n+\beta_n\mu}{\gamma_n+\delta_n\mu}\Gm_n(\mu),
\qquad n\in\Zp.
\end{equation}
We call weight functions that possess property \R{1.4} {\em M\"obius
quotient functions (MQF).\/}

It has been proved in \cite{IN2} that, as long as $\omega$ is a
M\"obius quotient function, the transformation
\R{1.2} assumes the form
\begin{equation}\label{1.5}
\GT \left\{\sum_{k=0}^n f_k \prod_{j=0}^{k-1} g_j(x) \prod_{j=k}^{n-1}
h_j(x) \right\}=\sum_{k=0}^n f_k x^k,
\end{equation}
where
$$g_k(\mu)=\alpha_k+\mu\beta_k,\quad h_k(\mu)=\gamma_k+\mu
\delta_k,\qquad k\in\Zp,$$
and where $f_0,f_1,\ldots,f_n$ are arbitrary.

Several MQFs that share strict sign consistency have been presented in
\cite{IN2}, for example $\omega_1(x,\mu)=x^\mu$,
$\omega_2(x,\mu)=\mu^x$ and $\omega_2(x,\mu)=x^{\log_q\mu}$,
$q\in(0,1)$. Moreover, the theme of \cite{IN5} was the identification
of all possible weight functions $\alpha(x)$ such that
$\omega_j(x,\mu)\alpha(x)$ is MQF and SSC for some $j\in\{1,2,3\}$.
The present paper takes yet another step toward our ultimate goal of
identifying all MQFs.

We commence our analysis in Section 2 with the assumption that the
coefficients $\alpha_n$, $\beta_n$, $\gamma_n$ and $\delta_n$ are
linear in $n$. This yields a linear ordinary differential equation
with polynomial coefficients that must be satisfied by $\omega$ and
that can be solved explicitly. We demonstrate that, up to a
transformation of variables, nothing can be gained beyond a known
transformation from \cite{IN2}.

In Section 3 we expand our framework to polynomial coefficients
$\alpha_n$, $\beta_n$, $\gamma_n$, $\delta_n$. The weight function
$\omega$ must again obey a linear ordinary differential equation,
which we solve explicitly in terms of a generalized hypergeometric
function.

Finally, in Section 4 we return to the representation \R{1.5}, in
order to derive an alternative form of the transformation
and of the underlying biorthogonal polynomial. Let
\begin{equation}\label{1.6}
\tilde{q}(x)=\frac{\prod_{k=1}^n (x-\mu_k)}{\prod_{k=0}^{n-1}
(\gamma_k +x\delta_k)}.
\end{equation}
We prove that
$$p_n(x;\mu_1,\mu_2,\ldots,\mu_n)=\sum_{k=0}^n
F_k[\tilde{q}_0,\tilde{q}_1,\ldots, \tilde{q}_k] x^k,$$
where $\tilde{q}_k=\tilde{q}(-\alpha_k/\beta_k)$, $k=0,1,\ldots,n$,
and $F_0,F_1,\ldots$ are a generalization of the familiar divided
differences. 

\section{Linear M\"obius quotient functions}

Let us suppose that the moment sequence of the weight function
$\omega$ obeys \R{1.4} with
\begin{equation}\label{2.1}
\alpha_n=a_0-n a_1,\quad \beta_n=b_0-n b_1,\quad \gamma_n=c_0-n
c_1,\quad \delta_n=d_0-n d_1,\qquad n\in\Zp.
\end{equation}
We further assume that $\omega(\,\cdot\,,\mu)\in C^1[0,\infty)$ and
that $\omega(0,\mu)\equiv0$ -- the last assumption does not lead to
severe loss of generality, since it can be always forced by shifting
$x\mapsto x+\varepsilon$ for some $\varepsilon>0$, smoothly
continuing $\omega$ for $x\in[0,\varepsilon)$ so that it vanishes at
the origin and ultimately driving $\varepsilon$ to zero. We say that
$\omega$ is a {\em linear\/} MQF.\footnote{The phrase `linear' refers,
of course, to the coefficients $\alpha_n$, $\beta_n$ etc., not to
$\omega$ itself!}

We rewrite \R{1.4} in the form
\begin{equation}\label{2.2}
(\gamma_n+\beta_n\mu)\int_0^\infty x^{n+1}\omega(x,\mu)\D
x=(\alpha_n+\beta_n \mu) \int_0^\infty x^n \omega(x,\mu)\D x,\qquad
\mu_-<\mu<\mu_+,
\end{equation}
multiply by $t^n$ and sum for $n=0,1,\ldots$. The outcome is the
identity
\begin{eqnarray}
&&\int_0^\infty \{ (a_0-c_0x)+\mu(b_0-d_0x)\}\omega(x,\mu)\E^{-tx}\D
x \label{2.3}\\
&=&-t\int_0^\infty \{(a_1-c_1x)+\mu(b_1-d_1x)\}x\omega(x,\mu) \E^{-tx}
\D x. \nonumber
\end{eqnarray}
Let
\begin{equation}\label{2.4}
\omega_\ell(x,\mu)=x^\ell \{(a_\ell-c_\ell x)+\mu(b_\ell-d_\ell x)\}
\omega(x,\mu),\qquad \ell=0,1,\quad \mu\in(\mu_-,\mu_+).
\end{equation}
We can thus rewrite \R{2.2} as
$${\cal L}\{\omega_0(\,\cdot\,,\mu)\}(t)=-t{\cal
L}\{\omega_1(\,\cdot\,,\mu)\}(t),$$
where $\cal L$ is the Laplace transform. However, for every $f\in
C^1$,
$${\cal L}\{f'\}(t)=t{\cal L}\{f\}(t)-f(0)$$
and we use the smoothness of $\omega$ and $\omega_1(0,\mu)\equiv0$ to
argue that 
$${\cal L}\left\{\frac{\partial}{\partial x} \omega_1(\,\cdot\,,\mu)
\right\}(t) ={\cal L}\{\omega_0(\,\cdot\,,\mu)\}(t).$$
Taking inverse transforms, we thus deduce that
\begin{equation}\label{2.5}
\frac{\partial}{\partial x}\omega_1(x,\mu)=-\omega_0(x,\mu),\qquad
x\geq0,\quad \mu\in(\mu_+,\mu_-).
\end{equation}
We finally substitute the explicit values of $\omega_0$ and $\omega_1$
from \R{2.3}, and this yields a linear ordinary differential equation
for $\omega$,
\begin{eqnarray}
&&x\{(a_1+\mu b_1)-(c_1+\mu d_1)x\}\frac{\partial}{\partial x}
\omega(x,\mu) \label{2.6}\\
&+&\left\{ [(a_0+\mu b_0)+(a_1+\mu b_1)]-[(c_0+\mu
d_0)+2(c_1+\mu d_1)]x\right\}\omega(x,\mu). \nonumber
\end{eqnarray}
This can be solved and, in tandem with the initial condition
$\omega(0,\mu)\equiv0$, results in a unique value of $\omega$.

For example, letting
$$\alpha_n=n,\quad \beta_n\equiv1,\quad \gamma_n=n+1,\quad
\delta_n\equiv1,\qquad n\in\Zp,$$
we obtain
$$a_0=0,\; a_1=-1,\quad b_0=1,\;b_1=0,\quad c_0=1,\;c_1=-1,\quad
d_0=1,\;d_1=0$$
and \R{2.6} becomes
$$\frac{\partial}{\partial
x}\omega(x,\mu)=\frac{\mu-1}{x}\omega(x,\mu),$$
with the solution (up to a multiplicative constant)
$\omega(x,\mu)=x^{\mu-1}$ (where $x>0$, $\mu>-1$) -- the Jacobi
transformation (with $\alpha=0$) from \cite{IN2}.

With greater generality, we observe that, up to a multiplicative
constant, the nonzero solution of
$$x(\sigma_0-\sigma_1x)y'+(\rho_0-\rho_1x)y=0,\qquad y(0)=0,$$
is
$$y(x)=x^{-\rho_0/\sigma_0} (\sigma_0-\sigma_1 x)^{\rho_0/\sigma_0
-\rho_1/\sigma_1}$$
-- subject, of course, to the above expression being well defined,
i.e.\ to $\rho_0/\sigma_0<0$ and either
$$\frac{\rho_0}{\sigma_0}\geq\frac{\rho_1}{\sigma_1},\qquad \sigma_0\geq
\sigma_1 x$$
or 
$$\frac{\rho_0}{\sigma_0}\leq\frac{\rho_1}{\sigma_1},\qquad \sigma_0<
\sigma_1 x$$

In our case 
$$\frac{\rho_0}{\sigma_0}=1+{\frac{a_0+\mu b_0}{a_1+\mu b_1}},$$
hence $\rho_0/\sigma_0<0$ is equivalent to
$$\frac{a_0+\mu b_0}{a_1+\mu b_1}+1<0.$$
Moreover,
$$\frac{\rho_0}{\sigma_0}-\frac{\rho_1}{\sigma_1}=-1+\frac{a_0+\mu
b_0}{a_1+\mu b_1}-\frac{c_0+\mu d_0}{c_1+\mu d_1}.$$

Particular substitutions readily yield an explicit form of $\omega$.
For example, we can use 
$${\bf a}=\left[\begin{array}{r}1\\0\end{array}\right],\quad {\bf
b}=\left[\begin{array}{r}0\\-1 \end{array}\right],\quad {\bf
c}=\left[\begin{array}{r}1\\0 \end{array}\right],\quad {\bf d}=\left[
\begin{array}{r}1\\-1\end{array}\right],$$
and this results in
$$\omega(x,\mu)=x^{1-\frac{1}{\mu}},\qquad \mu>1.$$

More general is a choice of ${\bf a},{\bf b},{\bf c}$ and ${\bf d}$
so that $\omega(x,\mu)=cx^{\nu(\mu)}$ for some $c>0$ and a function
$\nu(\mu)\not\equiv{\rm const}$. Thus, we need
$\rho_0\sigma_1=\rho_1\sigma_0$ and this, after some manipulation,
results in 
\begin{eqnarray*}
a_1c_0&=&(a_0-a_1)c_1,\\
b_1d_0&=&(b_0-b_1)d_1,\\
(a_0b_1-a_1b_0)(a_1d_1-b_1c_1)&=&0.
\end{eqnarray*}
Suppose first that $a_0b_1=a_1b_0$. Then
$$\frac{\rho_0}{\sigma_0}\equiv1+\frac{a_0}{a_1},$$
hence $\nu\equiv{\rm const}$, a contradiction. We thus deduce that
$a_1d_1=b_1c_1$. By letting $c_1=\kappa a_1$ we deduce the explicit expression
\begin{equation}\label{2.7}
d_1=b_1\kappa, \quad
c_0=(a_0-a_1)\kappa,\quad 
d_0=(b_0-b_1)\kappa,
\end{equation}
where $\kappa\neq0$ is arbitrary -- as are $a_0,a_1,b_0$ and $b_1$. We
deduce from \R{2.7} that
\begin{eqnarray*}
\alpha_n&=&a_0-a_1n,\\
\beta_n&=&b_0-b_1n,\\
\gamma_n&=&\kappa(a_0-a_1) -\kappa a_1 n,\\
\delta_n&=&\kappa(b_0-b_1)-\kappa b_1 n,\qquad n\in\Zp,
\end{eqnarray*}
and
$$\frac{\Gm_{n+1}(\mu)}{\Gm_n(\mu)}=\kappa^{-1} \frac{(a_0-a_1 n)
+(b_0-b_1 n)\mu}{((a_0-a_1)-a_1 n) +((b_0-b_1)-b_1n)\mu},$$
\begin{equation}\label{2.8}
\nu(\mu)=-1-\frac{a_0+\mu b_0}{a_1+\mu b_1}.
\end{equation}

We deduce that
\begin{eqnarray*}
g_n(\mu)&=&-(a_1+b_1\mu)\left(-\frac{a_0+b_0\mu}{a_1+b_1\mu}+n\right),\\
h_n(\mu)&=&-\kappa(a_1+b_1\mu)\left(-\frac{a_0+b_0\mu}
{a_1+b_1\mu}+n+1\right),\qquad n\in\Zp.
\end{eqnarray*}
Moreover, $a_1b_0\neq a_0 b_1$ readily implies \R{1.3}.

The function $\nu(\mu)$ being strictly monotone, we deduce strict sign
consistency of $\omega$ from strict sign consistency of the function
$x^\mu$ \cite{KS,PS}. Moreover,
$$\prod_{j=0}^{k-1}g_j(x)\prod_{j=k}^{n-1} h_j(x)
=(-1)^{n+1}\kappa^{n-k} (a_1+b_1x)^n \frac{
\left( -\frac{a_0+b_0 x}{a_1+b_1 x}\right)_{n+1}} {k-\frac{a_0+b_0
x}{a_1+b_1 x}},\qquad k=0,1,\ldots,n,$$
where
$$(a)_n=a(a+1)(a+2)\cdots(a+n-1),\qquad a\in\CC,\quad n\in\Zp,$$
is the familiar {\em Pocchammer symbol\/}. Substitution in \R{1.5}
thus yields the transformation
$$\GT\left\{\left(\prod_{\ell=0}^n ((a_1+b_1 x)\ell -(a_0+b_0
x))\right) \sum_{k=0}^n \frac{\kappa^{-k} f_k}{(a_1+b_1x)k-(a_0+b_0 
x)}\right\} =\sum_{k=0}^n f_k x^k.$$
This, however, can be converted, by a change of variable, to the
jacobi transformation with $\alpha=0$ \cite{IN2}, which can be
rendered in the form
$$\GT \left\{ (x)_{n+1} \sum_{k=0}^n
\frac{f_k}{k+x}\right\}=\sum_{k=0}^n f_k x^k.$$

The case $\rho_0\sigma_1\neq \rho_1\sigma_0$ can be resolved
explicitly by similar techniques. 

\section{Polynomial M\"obius quotient functions}

The approach of the last section can be generalized to cater for a
more extensive family of MQFs. Thus, suppose that the coefficients can
be expanded in the form
\begin{eqnarray*}
\alpha_n&=&\sum_{\ell=0}^\infty a_\ell (-n)_\ell,\\
\beta_n&=&\sum_{\ell=0}^\infty b_\ell (-n)_\ell,\\
\gamma_n&=&\sum_{\ell=0}^\infty c_\ell (-n)_\ell,\\
\delta_n&=&\sum_{\ell=0}^\infty d_\ell (-n)_\ell,
\end{eqnarray*}
and that the functions
$$a(z)=\sum_{\ell=0}^\infty a_\ell z^\ell, \quad
b(z)=\sum_{\ell=0}^\infty b_\ell z^\ell, \quad
c(z)=\sum_{\ell=0}^\infty c_\ell z^\ell, \quad
d(z)=\sum_{\ell=0}^\infty d_\ell z^\ell$$
are analytic about the point $z=0$. We further assume that
$\omega(\,\cdot\,,\mu)\in C^\infty [0,\infty)$ and that
\begin{equation}\label{3.1}
\frac{\partial^\ell}{\partial x^\ell}\omega(0,\mu)\equiv0,\qquad
\ell=0,1,\ldots.
\end{equation}

The last condition, in tandem with infinite differentiability, is
restrictive. This can be alleviated along the lines of the last
section, by letting
$$\omega_\varepsilon(x,\mu)=\left\{ \begin{array}{lcl}
\omega(x-\varepsilon,\mu) & \qquad & :\varepsilon\leq x,\\
\sum_{\ell=0}^\infty \theta_\ell \E^{-(\ell+1)x^{-2}} & \qquad &
:0< x\leq \varepsilon,    \end{array}\right.$$
where $\theta_0,\theta_1,\ldots$ are chosen so as to ensure infinite
differentiability at $x=\varepsilon$,\footnote{The existence of such
coefficients follows easily from the M\"untz theorem \cite{Da}.} and letting
$\varepsilon\rightarrow0$ by the end of the procedure. The situation is
simplified if -- as we do further in this section -- it is assumed
that $a,b,c$ and $d$ are polynomials, since there it is enough to
insist that $\omega(\,\cdot\,,\mu)\in C^r[0,\infty)$ and that \R{3.1}
holds for $\ell=0,1,\ldots,s-1$, where
$$s=\max\{\deg a,\deg b,\deg c,\deg d\}.$$

We again multiply \R{2.2} by $t^n$ and sum up for all integer
nonnegative $n$. This results in the identity
\begin{equation}\label{3.2}
\int_0^\infty x\{c(tx)+\mu d(tx)\}\omega(x,\mu)\E^{-tx}\D
x=\int_0^\infty \{a(tx)+\mu b(tx) \omega(x,\mu) \E^{-tx}\D x.
\end{equation}
Expanding
$$x\{c(tx)+\mu d(tx)\} \qquad \mbox{and} \qquad \{a(tx)+\mu b(tx)\}$$
in powers of $t$, we deduce from \R{3.2} that
\begin{equation}\label{3.3}
\sum_{\ell=0}^\infty t^\ell \int_0^\infty \omega_\ell(x,\mu)\E^{-tx}\D
x=0, 
\end{equation}
where $\omega_\ell$ has been given by \R{2.4}. Since \R{2.4} and
\R{3.1} combine to produce
$$\frac{\partial^j}{\partial x^j}\omega_\ell(0,\mu)=0,\qquad
j=0,1,\ldots,\ell-1,$$
we conclude that
$$t^\ell\int_0^\infty \omega_\ell(x,\mu)\E^{-tx}\D x={\cal L}\left\{
\frac{\partial^\ell} {\partial x^\ell}
\omega_\ell(x,\mu)\right\}(t),\qquad \ell\in\Zp.$$
Taking inverse transforms, \R{3.3} becomes
\begin{equation}\label{3.4}
\sum_{\ell=0}^\infty \frac{\partial^\ell}{\partial x^\ell}
\omega_\ell(x,\mu)=0.
\end{equation}

Let us set
$$P_\ell(\mu)=a_\ell+b_\ell\mu,\quad
Q_\ell(\mu)=c_\ell+d_\ell\mu,\qquad \ell\in\Zp,$$
therefore
$$\omega_\ell(x,\mu)=(P_\ell(\mu)-Q_\ell(\mu)x)x^\ell
\omega(x,\mu),\qquad \ell\in\Zp.$$
Substitution into \R{3.4} and differentiation with the Leibnitz rule
result in
\begin{equation}\label{3.5}
\sum_{j=0}^\infty \frac{1}{j!(j+1)!} \left((j+1) \sum_{\ell=0}^\infty
\frac{P_{\ell+j}(\mu)}{\ell!} -x\sum_{\ell=0}^\infty
(\ell+j+1)\frac{Q_{\ell+j}(\mu)}{\ell!} \right) x^j
\frac{\partial^j}{\partial x^j} \omega(x,\mu)=0.
\end{equation}
Note that both
$$\Gp_j(\mu):=\sum_{\ell=0}^\infty
\frac{P_{\ell+j}(\mu)}{\ell!}=\sum_{\ell=0}^\infty
\frac{a_{\ell+j}}{\ell!} +\mu \sum_{\ell=0}^\infty
\frac{b_{\ell+j}}{\ell!}$$
and
$$\Gq_j(\mu):=\sum_{\ell=0}^\infty (\ell+j+1)
\frac{Q_{\ell+j}(\mu)}{\ell!}=\sum_{\ell=0}^\infty (\ell+j+1)
\frac{c_{\ell+j}}{\ell!} +\mu \sum_{\ell=0}^\infty (\ell+j+1)
\frac{d_{\ell+j}}{\ell!}$$
are linear in $\mu$.

We henceforth assume that $a,b,c,d\in{\Bbb P}_s[\mu]$, designating the
underlying $\omega$ as {\em polynomial\/} MQF. Therefore \R{3.5} is an
ordinary differential equation of degree $s$, equipped with zero
initial conditions at the origin.

The equation \R{3.5} is a
special case of the $s$-degree ODE
\begin{equation}\label{3.6}
\sum_{\ell=0}^s (p_\ell-q_\ell x)x^\ell y^{(\ell)}=0.
\end{equation}
We express the solution of \R{3.6} in the form
$$y(x)=x^\theta \sum_{n=0}^\infty y_n x^n.$$
Substitution into \R{3.6} and elementary manipulation affirm that, for
arbitrary $y_0$, $\theta$ must be a zero of
\begin{equation}\label{3.7}
\sum_{\ell=0}^s (-1)^\ell p_\ell (-\theta)_\ell=0
\end{equation}
(thus, provided that $p_s\neq0$, there are $s$ solutions in $\CC$, the
right number!), whereas the coefficients $y_n$ are obtained from the
recurrence
\begin{equation}\label{3.8}
y_n=\frac{\sum_{\ell=0}^s (-1)^\ell q_\ell (-n+1-\theta)_\ell}
{\sum_{\ell=0}^s (-1)^\ell p_\ell (-n-\theta)_\ell}y_{n-1},\qquad
n=1,2,\ldots.
\end{equation}
Let
$$P(x):=\sum_{\ell=0}^s(-1)^\ell p_\ell (-x-\theta)_\ell,\qquad
Q(x):=\sum_{\ell=0}^s (-1)^\ell q_\ell (-\theta-x)_\ell.$$
Therefore, \R{3.7} is equivalent to $P(0)=0$, whereas \R{3.8} converts
to
$$y_n=\frac{Q(n-1)}{P(n)}y_{n-1}=y_0\prod_{\ell=1}^n
\frac{Q(\ell-1)}{P(\ell)}, \qquad n=1,2,\ldots.$$
Suppose that
$$P(x)=p_*\prod_{j=1}^s (x-\eta_j),\qquad Q(x)=q_* \prod_{j=1}^s
(x-\zeta_j),$$
where, without loss of generality, $\eta_s=0$ (recall that $P(0)=0$).
Then 
$$\prod_{\ell=1}^n P(\ell)=p_*^n \prod_{j=1}^s (1-\eta_j)_n,\quad
\prod_{\ell=0}^{n-1} Q(\ell)=q_*^n \prod_{j=1}^s (-\zeta_j)_n$$
therefore
$$y_n=\left(\frac{q_*}{p_*}\right)^n \frac{\prod_{j=1}^s (-\zeta_j)_n}
{\prod_{j=1}^s (1-\eta_j)_n} y_0,\qquad n\in\Zp.$$
Letting, without loss of generality, $y_0=1$, we deduce that
$$y(x)=x^\theta {}_sF_{s-1} \left[ \begin{array}{l}
-\zeta_1,-\zeta_2,\ldots,-\zeta_s;\\
1-\eta_1,1-\eta_2,\ldots,1-\eta_{s-1};\end{array} \frac{q_*}{p_*}x
\right],$$
a {\em generalized hypergeometric\/} function \cite{Ra}.

Of course, we need $\theta\geq s$, otherwise the requisite number of
derivatives at the origin fails to vanish -- note that we need only
$\ell=0,1,\ldots,s-1$ in \R{3.1}.

The last result can be somewhat generalized by observing that the
degree of either $P$ or $Q$ may be strictly less than $s$ and
repeating our analysis without any significant changes.

\vspace{8pt}
\noindent {\bf Theorem 1} Let $\alpha_n,\beta_n,\gamma_n,\delta_n\in
{\Bbb P}_s[n]$, $n\in\Zp$, and suppose further that
$\omega(\,\cdot\,,\mu)\in C^s[0,\infty)$ and that $\partial^\ell
\omega(0,\mu)/\partial x^\ell=0$, $\ell=0,1,\ldots,s-1$. Then there
exist nonnegative integers $s_1,s_2\leq s$ and real numbers
$\theta$, $\zeta_1,\zeta_2,\ldots,\zeta_{s_1}$, $\eta_1,\eta_2,
\ldots,\eta_{s_2-1}$, $\nu$, all functions of $\mu$, such that
\begin{equation}\label{3.9}
\omega(x,\mu)=x^\theta {}_{s_1}F_{{s_2}-1} \left[ \begin{array}{l}
-\zeta_1,-\zeta_2,\ldots,-\zeta_{s_1};\\
1-\eta_1,1-\eta_2,\ldots,1-\eta_{{s_2}-1};\end{array} \nu x \right].
\end{equation}
$\Box$

\vspace{8pt}
We remark that, needless to say, not every ${}_{s_1} F_{{s_2}-1}$ function
will do in \R{3.9}, since we require $\omega(x,\mu)\geq0$ for all
$x\in[0,\infty)$, $\mu_-<\mu<\mu_+$. Moreover, the more interesting
case is when, in addition, $\omega$ is SSC.

As an illustration of our technique, we consider in detail the case
when the quotient \R{1.4} is linear in $\mu$, i.e.\ $\gamma_n\equiv1$,
$\delta_n\equiv0$, and assume in addition that $s=2$. Therefore
$$\alpha_n=a_0-a_1n+a_2n(n-1),\quad \beta_n=b_0-b_1n+b_2n(n-1),\qquad
n\in\Zp,$$
and we have
$$\begin{array}{rclcrcl}
P_0(\mu) &=& a_0+b_0\mu, & \qquad & Q_0(\mu) &\equiv& 1,\\
P_1(\mu) &=& a_1+b_1\mu, & \qquad & Q_1(\mu) &\equiv& 0,\\
P_2(\mu) &=& a_2+b_2\mu, & \qquad & Q_2(\mu) & \equiv& 0.
  \end{array}$$
In other words,
$$\begin{array}{rclcrcl}
\Gp_0(\mu) &=&\left(a_0+a_1+\Frac12 a_2\right)+\left(b_0+b_1+\Frac12 b_2
\right)\mu, & \qquad & \Gq_0(\mu) &\equiv& 1,\\
\Gp_1(\mu) &=& (a_1+a_2)+(b_1+b_2)\mu, &\qquad& \Gq_1(\mu) &\equiv&
0,\\
\Gp_2(\mu) &=& a_2+b_2\mu, &\qquad& \Gq_2(\mu) &\equiv& 0
  \end{array}$$
and the differential equation becomes
\begin{eqnarray*}
\Frac{1}{4} (a_2+b_2\mu)x^2\frac{\partial^2\omega(x,\mu)}{\partial
x^2}+((a_1+a_2)+(b_1+b_2)\mu) x\frac{\partial\omega(x,\mu)}{\partial
x}&&\\
\mbox{}+\left(\left(a_0+a_1+\Frac12 a_2\right)+\left(b_0+b_1+\Frac12 b_2
\right)\mu-x\right)\omega(x,\mu)&=&0,
\end{eqnarray*}
with zero initial conditions at $x=0$. Since
$$P(x)=\Gp_0(\mu)+\Gp_1(\mu)(x+\theta)+\Frac14\Gp_2(\mu)(x+\theta)
(x+\theta-1)$$
and we require $P(0)=0$, we choose $\theta=\theta(\mu)$ as a root of
the quadratic
\begin{equation}\label{3.10}
\Frac14(a_2+b_2\mu)\theta^2+\left(\left(a_1+\Frac34 a_2\right)+\left(
b_1+\Frac34 b_2\right)\mu\right)\theta+\left(a_0+a_1+\Frac12
a_2\right)+\left(b_0+b_1+\Frac12 b_2 \right)\mu=0.
\end{equation}
Therefore
$$\eta_1=1-2\theta-4\frac{\Gp_1}{\Gp_2}$$
and $s_2=2$. On the other hand, $Q(x)\equiv1$, therefore $s_1=0$.
Moreover, $\nu=4/\Gp_2$ and we deduce from Theorem 1 that
\begin{equation}\label{3.11}
\omega(x,\mu)=x^\theta {}_0F_1 \left[ \begin{array}{l}
\mbox{---};\\2\theta+4\Gpp_1/\Gpp_2; \end{array} \frac{4x}{\Gp_2}
\right].
\end{equation}
Recall that a {\em modified Bessel function\/} can be represented in
the form \cite{Ra}
$$I_\tau(z)=\frac{\left(\frac12 z\right)^\tau}{\Gamma(\tau+1)} {}_0F_1
\left[ \begin{array}{l}\mbox{---};\\\tau+1; \end{array} \Frac14
z^2\right].$$
Let us assume that $\mu_-,a_2,b_2\geq0$. Then, dispensing with a
multiplicative factor, we have
$$\omega(x,\mu)=x^{\frac12-2\Gpp_1/\Gpp_2}
I_{2\theta+4\Gpp_1/\Gpp_2-1}\left(4\left(\frac{x}{\Gp_2}\right)^{\frac12}
\right).$$
As long as $\theta$ is real and $\geq1$ we are assured that the
initial conditions are satisfied and that $\omega(x,\mu)\geq0$ for all
$x\geq0$. 

Let us consider the special choice
$$a_0=a_1=b_0=b_2=0,\qquad a_2>0,\qquad b_1\neq0.$$
In other words, $\alpha_n=a_2n(n-1)$ and $\beta_n=-b_1n$ for $n\in\Zp$.
Letting $\lambda=-b_1\mu/a_2$, \R{3.10} becomes
$$\theta^2+(3-4\lambda)\theta+(2-4\lambda)=0,$$
with the solutions $\theta=-1$ and $\theta=4\lambda-2$. Since we
require $\theta\geq1$, we choose the second solution and impose
$\lambda\geq\frac34$. The outcome is
$$\omega(x,\mu)=x^{2\lambda-\frac32} I_{4\lambda-1}
\left(4\left(\frac{x}{a_2}\right)^{\frac12}\right),$$
where $\mu$ need be restricted consistently with $b_1\mu+\frac34
a_2\leq0$. This can be simplified, e.g.\ by considering the variables
$$\tilde{x}=4\left(\frac{x}{a_2}\right)^{\frac12},\qquad
\tilde{\mu}=\lambda-\Frac34,$$ 
instead of $x$ and $\mu$ respectively. This results in
$$\omega(\tilde{x},\tilde{\mu})=\tilde{x}^{\tilde{\mu}}
I_{\tilde{\mu} +1}(\tilde{x}),\qquad \tilde{x}\geq0,\quad
\tilde{\mu}>0.$$

\section{An inversion formula for a biorthogonal representation}

The main result underpinning the transformation \R{1.5} is Lemma 5 of
\cite{IN2}, namely that, as long as 
\begin{equation}\label{4.1}
\prod_{k=1}^n (x-\mu_k)=\sum_{k=0}^n f_k \prod_{j=0}^{k-1} (\alpha_j
+x\beta_j) \prod_{j=k}^{n-1} (\gamma_j +x\delta_j),
\end{equation}
and the moments of $\omega$ $\{\Gm_s(\mu)\}_{s\in\Zp}$ are consistent
with \R{1.4}, the $n$th biorthogonal polynomial is
\begin{equation}\label{4.2}
p_n(x;\mu_1,\mu_2,\ldots,\mu_n)=\sum_{k=0}^n f_k x^k.
\end{equation}

Our purpose in the present section is to express the coefficients
$f_0,f_1,\ldots,f_n$ 
explicitly in terms of the parameters $\mu_1,\mu_2,\ldots,\mu_n$. Let
$$\tilde{q}(x)=\frac{\prod_{k=1}^n (x-\mu_k)}{\prod_{k=0}^{n-1}
(\gamma_k +x\delta_k)}.$$
Then \R{4.1} yields
\begin{equation}\label{4.3}
\tilde{q}(x)=\sum_{k=0}^n f_k \prod_{j=0}^{k-1} \Gm_j(x).
\end{equation}
We let $\lambda_j=-\alpha_j/\beta_j$, $j=0,1,\ldots,n$ and observe
that, by \R{1.4}, $\Gm_j(\lambda_\ell)=0$ for $\ell=0,1,\ldots,j-1$.
Therefore, letting $x=\lambda_\ell$ in \R{4.3}, we obtain
\begin{equation}\label{4.4}
\tilde{q}_\ell=\sum_{k=0}^\ell f_k \prod_{j=0}^{k-1}
\Gm_j(\lambda_\ell),\qquad \ell=0,1,\ldots,n,
\end{equation}
where $\tilde{q}_\ell=\tilde{q}(\lambda_\ell)$, $\ell=0,1,\ldots,n$.

We assume that $\Gm_j(\lambda_\ell)\neq0$ for $\ell> j$. Letting
$\ell=0,1,\ldots$ in \R{4.4}, we obtain 
\begin{eqnarray*}
f_0&=&\tilde{q}_0,\\
f_1&=&\frac{\tilde{q}_1-\tilde{q}_0}{\Gm_0(\lambda_1)},\\
f_2&=&\frac{\displaystyle \frac{\tilde{q}_2-\tilde{q}_0}{\Gm_0(\lambda_2)}
-\frac{\tilde{q}_1-\tilde{q}_0}{\Gm_0(\lambda_1)}}{\Gm_1(\lambda_2)}
\end{eqnarray*}
etc. To obtain the general form of $f_k$, we need to generalize the
concept of divided differences. Given a sequence
$\{r_\ell\}_{\ell\in\Zp}$, we define
\begin{eqnarray*}
F_0[r_\ell]&=&r_\ell,\\
F_k[r_\ell,r_{\ell+1},\ldots,r_{\ell+k}]&=&\frac{F_{k-1}[r_\ell,\ldots,
r_{\ell+k-2},r_{\ell+k}-F_{k-1}[r_\ell,\ldots,r_{\ell+k-1}]}{\Gm_{k-1}
(\lambda_{\ell+k})}, \quad k=1,2,\ldots.
\end{eqnarray*}

\vspace{8pt}
\noindent {\bf Lemma 2}
$f_k=F_k[\tilde{q}_0,\tilde{q}_1,\ldots,\tilde{q}_k]$,
$k=0,1,\ldots,n$.

\vspace{6pt}
{\em Proof.\/} By induction on $\ell$. The statement is certainly true
for $\ell=0$ and we assume it for $\ell-1$. Therefore, it follows from
\R{4.4} that
$$\tilde{q}_\ell=\tilde{q}_0+\sum_{k=1}^{\ell-1}
F_k[\tilde{q}_0,\ldots,\tilde{q}_k] \prod_{j=0}^{k-1}\Gm(\lambda_j)
+d_\ell \prod_{j=0}^{\ell-1} \Gm_j(\lambda_\ell).$$
Therefore,
\begin{eqnarray*}
f_\ell&=&\frac{\tilde{q}_\ell-\tilde{q}_0}{\prod_{j=0}^{\ell-1}
\Gm_j(\lambda_\ell)}-\sum_{k=1}^{\ell-1}
\frac{F_k[\tilde{q}_0\ldots,\tilde{q}_k]} {\prod_{j=k}^{\ell-1}
\Gm_j(\lambda_\ell)}\\
&=&\frac{F_1[\tilde{q}_0,\tilde{q}_\ell]}{\prod_{j=1}^{\ell-1}
\Gm_j(\lambda_\ell)} -\sum_{k=1}^{\ell-1}
\frac{F_k[\tilde{q}_0\ldots,\tilde{q}_k]} {\prod_{j=k}^{\ell-1}
\Gm_j(\lambda_\ell)}.
\end{eqnarray*}
We will prove by induction on $s=1,2,\ldots,\ell-2$ that
\begin{equation}\label{4.5}
f_\ell=\frac{F_s[\tilde{q}_0,\ldots,\tilde{q}_{s-1},\tilde{q}_\ell]}
{\prod_{j=s}^{\ell-1} \Gm_j(\lambda_\ell)} -\sum_{k=s}^{\ell-1}
\frac{F_k[\tilde{q}_0,\tilde{q}_1,\ldots,\tilde{q}_k]}
{\prod_{j=k}^{\ell-1} \Gm_j(\lambda_\ell)}.
\end{equation}
This is certainly true for $s=1$ and, assuming that \R{4.5} is valid
for $s\geq1$, we have
\begin{eqnarray*}
f_\ell&=&\frac{F_s[\tilde{q}_0,\ldots,\tilde{q}_{s-1},\tilde{q}_\ell]
-F_s[\tilde{q}_0,\ldots,\tilde{q}_{s-1},\tilde{q}_s]}{\prod_{j=s}^{\ell-1}
\Gm_j(\lambda_\ell)} -\sum_{k=s+1}^{\ell-1} \frac{F_k[\tilde{q}_0,
\tilde{q}_1,\ldots,\tilde{q}_k]}{\prod_{j=k}^{\ell-1}
\Gm_j(\lambda_j)}\\
&=&\frac{F_{s+1}[\tilde{q}_0,\ldots,\tilde{q}_s,\tilde{q}_{s+1}]}
{\prod_{j=s+1}^{\ell-1} \Gm_j(\lambda_\ell)} -\sum_{k=s+1}^{\ell-1}
\frac{F_k[\tilde{q}_0,\tilde{q}_1,\ldots,\tilde{q}_k]}
{\prod_{j=k}^{\ell-1} \Gm_j(\lambda_j)}.
\end{eqnarray*}
This proves that \R{4.5} is true. We now let $s=\ell-2$, therefore
\begin{eqnarray*}
f_\ell&=&\frac{F_{\ell-1}[\tilde{q}_0,\ldots,\tilde{q}_{\ell-2},
\tilde{q}_\ell]-F_{\ell-1}[\tilde{q}_0,\ldots,\tilde{q}_{\ell-2},
\tilde{q}_{\ell-1}]}{\Gm_{\ell-1}(\lambda_\ell)}\\
&=&F_\ell[\tilde{q}_0,\tilde{q}_1,\ldots,\tilde{q}_\ell].
\end{eqnarray*}
This completes the proof. \QED

\vspace{8pt}
\noindent {\bf Theorem 3} Let the moments of $\omega$ satisfy \R{1.4}
and assume that $\beta_\ell\neq0$, $\ell\in\Zp$, and
$\alpha_j\beta_\ell\neq \alpha_\ell\beta_j$ for all $\ell,j\in\Zp$,
$\ell\neq j$. Set
$$\tilde{q}_\ell=\frac{\prod_{k=1}^n (\alpha_\ell+\beta_\ell\mu_k)}
{\prod_{k=0}^{n-1} (\alpha_\ell \delta_k-\beta_\ell \gamma_k)},\qquad
\ell=0,1,\ldots,n.$$
Then
\begin{equation}\label{4.6}
p_n(x;\mu_1,\mu_2,\ldots,\mu_n)=\sum_{k=0}^n
F_k[\tilde{q}_0,\tilde{q}_1,\ldots,\tilde{q}_k] x^k.
\end{equation}

\vspace{6pt}
{\em Proof.\/} Follows at once from the lemma, since the above
definition of $\{\tilde{q}_\ell\}_{\ell=0}^n$ is consistent with
$\tilde{q}_\ell =\tilde{q}(\lambda_\ell)$. \QED

\end{document}